\documentclass{amsart}
\usepackage{amssymb,amsmath,theorem,euscript}

\usepackage[X2,T2A]{fontenc}
\usepackage[cp1251]{inputenc}
\usepackage[russian, english]{babel}

\newcommand{\и}{{\fontencoding{X2}\selectfont\cyrii}} 

\newcommand{\е}{{\fontencoding{X2}\selectfont\cyryat}} 

\newcounter{sec}

\newcounter{punct}[sec]

\def\punct{\refstepcounter{punct}{\arabic{sec}.%
\arabic{punct}.  }}

\def\punct{\refstepcounter{punct}{\arabic{sec}.\arabic{punct}.  }}

\def\COUNTERS{\addtocounter{sec}{1}
              \setcounter{punct}{0}
          \setcounter{equation}{0}
          }
          
          \def\sm{\smallskip}

\def\R{\mathbb{R}}
\def\Z{\mathbb{Z}}

\def\cS{\mathcal S}

\def\wt{\widetilde}

\def\ge{\geqslant}
\def\le{\leqslant}
\renewcommand\phi{\varphi}

\renewcommand{\Re}{\mathop {\mathrm {Re}}\nolimits}

\begin{document}
	
	
\begin{center}
	
	\large \bf

The double gamma function and Vladimir Alekseevsky

\bigskip

\sc

Yury A. Neretin
	






\end{center}

\bigskip

{\small 
This  paper is about a forgotten  function and a forgotten mathematician.
The double gamma function  is now an important special function,
which appears for different reasons in many branches of mathematics
and in mathematical physics,
as it satisfies a large and amusing collection
of identities parallel to the classical gamma function. It was discovered and investigated in detail by Vladimir Alekseevsky in  1888-89.
We outline his starting point here:
 considering the Weierstrass product for the entire periodic function
   $\sin \pi x$ and  
 taking half of the factors corresponding to non-positive roots
 of the sine,
  we obtain the function $1/\Gamma(x)$.
Considering the Weierstrass product for the doubly quasiperiodic Jacobi theta function 
$\vartheta_1$
and taking the half of factors we come to the $q$-gamma function (which
appears, with a slight change in notation, in Eduard Heine's work). 
Considering the
quarter of the factors (corresponding to the positive quadrant
in the lattice $n_1\omega+n_2\omega_2$ of periods), 
Alekseevsky arrived at his double gamma function. 
 The work in this direction 
 was continued
by Ernest Barnes, Jean Beaupin,
Godfrey Hardy, and Vladimir Steklov in 1899--1907. After 1907,
publications on this topic stopped, and the double gamma appeared again 70 years later due to Takuro Shintani.
In this paper we discuss  Alekseevsky's seminal work and its genesis, the history
of the double gamma, and   Alekseevsky's  biography (1858-1916).


}


\bigskip

{\small



}

\bigskip

{\sc   \phantom{4000000000000000000000000.}Content.
\newline
1. Introduction. A forgotten function and a forgotten mathematician
\newline
2. Biography of Alekseevsky
\newline
3. Alekseevsky's work, 1888-89
\newline
4. On the history of the double gamma function
\newline 
5. Alekseevsky: historiography,  lost work, and publications
\newline
6. Conslusion.
\newline
\phantom{4.}References 
}

\section{Introduction. A forgotten function and a forgotten mathematician}

\COUNTERS

In 1996, Vladimir Fock spoke to me about the work \cite{Fad},  \cite{Fock0}, \cite{Kash}, \cite{Fock} being done on the quantization of Teichm\"uller spaces; this construction was essentially based on a couple amusing special functions called 'Barnes double sine' and 'Barnes double gamma'.
I read the original   papers of Ernest Barnes \cite{Bar1}--\cite{Bar5} and
observed that he referred to an earlier work
of Alexeiewsky. This name was unknown to me, and my immediate attempt 
to find an information about this person and his works  was unsuccessful. 

Many years later, collecting  material for a book
 about Nikolay Luzin,  I opened 
 Krulikovsky's '{\it History of the development of mathematics in Tomsk}'
 \cite{Kru}.
Here I was again met with the name Aleksevsky. I found his main work
and
 was very surprised reading it.
   I then  started a historical 
excavation, which  I now  present to the reader.  

\sm

The double gamma function   is now a well recognized special function
 used in different branches of mathematics and in mathematical physics.
  It is an entire function $G(x;\alpha)$
of a complex variable $x$ and  depending on a complex parameter $\alpha$,
which
satisfies the following functional relations:
\begin{align}
G(x+1;\alpha)&=\Gamma\Bigl(\frac x\alpha \Bigr)\, G(x;\alpha);
\label{eq:01}
\\
G(x+\alpha;\alpha)&=(2\pi)^{\frac{\alpha-1}{2}} \alpha^{-\frac{2x-1}{2}}
\,
\Gamma(x)\,G(x;\alpha),
\label{eq:02}
\end{align}
and is normalized by
$$
G(0;\alpha)=1.
$$
This function satisfies a big collection of properties parallel to
the classical 
gamma function.

A special case of $G(x,\alpha)$ corresponding to $\alpha=1$ was discovered by Hermann Kinkelin
in \cite{Kin}; notice that  in this case
the two equations \eqref{eq:01}, \eqref{eq:02} coincide. His remarkable work did not attract  interest and was forgotten for a long time. The double gamma for general $\alpha$ was introduced by Vladimir Alekseevsky%
   \footnote{Russian: Алексеевский; old Russian orthography: Алекс\еевск\ий. 
 I use the most natural English spelling of the name,
 which is used by his modern homonyms.
German and French transliterations used in 1884-1907
were Alexeiewsky,  Alexejewsky, Alex\'eievsky.} in a large memoir
\cite{Ale88}--\cite{Ale89}. This work was continued by 
Ernest Barnes \cite{Bar1}, \cite{Bar2}, \cite{Bar3}, \cite{Bar4},
\cite{Bar-cap}, \cite{Bar5}, Jean Beaupain \cite{Bea1903}, \cite{Bea1904}, \cite{Bea1907},
Godfrey  Hardy \cite{Har1}, \cite{Har2}, and
 Vladimir Steklov \cite{Ste}.  The chain of publications
 on the double gamma and the multiple gamma 
 was interrupted in 1907, and after which  they  were more or less
 forgotten for 70 years. 
 
 Alekseevsky (1858--1916), the discover of the double gamma,  was also
  forgotten.
 He was mentioned 
in a few Russian papers on the histories of mathematics
and  education, 
\cite{Gus0}, \cite{Gus}, \cite{Kru}, \cite{Bah}, \cite{Kru1},
\cite{Yush},
\cite{Gra2003}, 
\cite{Bel} (mainly, on the local history of mathematics in Kharkov and Tomsk). He studied and worked at Kharkov University (1877-1905)
and  Tomsk Technological Institute (1905-1916).
For the Kharkov period, we mainly have records of unsuccessful attempts to 
advance his career.   
In Tomsk he became a full professor of mechanics, and in 1907-1911 he was 
the director (rector) of the Institute. 

\sm

We briefly discuss his biography in Sect. 2. But our main topics
are his impressive central work and the difficulty of mathematical development.

In Sect.~3, we discuss  Alekseevsky's main work
\cite{Ale88}--\cite{Ale89} and its genesis. It was a non-obvious response
to Eduard Heine's investigations of $q$-hypergeometric functions
 (though, more precisely, 
Alekseevsky's starting point  was Paul Appell's work \cite{App} on
extending  the $q$-gamma).
 
In Sect. 4, we discuss the history of the double gamma. We start with
a brief exposition of Kinkelin's work. It was basically surpassed
in \cite{Ale88}--\cite{Ale89}, but it also contains important results, which were 
not repeated in Alekseevsky's memoir (the link between Kinkelin and Alekseevsky was later observed by Beaupin). We also try to understand the origin
of Kinkelin's work. It seems to have been an investigation in a 'speleological'
style, which is relatively common in theory of special functions. 

\sm

The forgetting of  important research%
\footnote{Opening the database MathSciNet, which is quite selective, I see 233 references to 
\cite{Bar1} starting 2002. At the very  least, this is unusual for a paper published
a hundred years before. Google Scholar shows 512 references to \cite{Bar5}.}  is in itself an interesting historical phenomenon. In this case, we see
impressive articles 
clearly written  in a common mathematical language; they were
available to  readers of the XX century as well as  readers of 1890-00s.
We are trying to understand this phenomenon.

 Sect. 5 contains three additions. Firstly, we briefly discuss historical sources on Alekseevsky and
historiography. 

Secondly, 
starting with \cite{Ale92}, \cite{Ale94}
 Alekseevsky  claimed many times over
his researches could be extended to a wider class of functions depending 
on several parameters. We analyze these claims --- clearly, 
he kept in mind
 the multiple gamma functions; definitely he did not publish  
theorems on this topic. However, these claims were known and made water run in circles.
Apparently, \cite{Bar5} was a response to a declaration in \cite{Ale94}.

 To be complete, I  present a list of his  textbooks on calculus
and mechanics published by universities (all research publications are listed in Bibliography,
\cite{Ale84-1}--\cite{Ale04}; I think that this list is complete).

\section{Biography of Alekseevsky}

\COUNTERS

 {\bf \punct Early biography.} Vladimir Petrovich
  Alekseevsky was born April 19 (Julian 07), 1858
  in  Novgorod  governorate%
  \footnote{Russian: `губерния' (`guberniya').}
 and baptized in the village Peret\"enka
  located near  Okulovka  station on
   the railway line Moscow--St.Petersburg.
 His father, Petr Osi\-po\-vich Alekseevsky, was a `praporshchik' [$\simeq$ warrant officer or sergeant major] of the Corps of tran\-sport engineers [Корпус инженеров путей сообщения], \cite{Sin}, \cite{Gra2003}.
  
 In 1877, Alekseevsky graduated from the Nizhny Novgorod  gymnasium. 
It was a good school, with many of its pupils going on to become well known, in particular,
 mathematicians Aleksandr Lyapunov and Ivan Privalov,
 physicist Grigory Landsberg%
  \footnote{
Landsberg in particular could have won the Nobel in 1930 (combinational scattering of light),
 but was not nominated for it.
 He also edited the famous high-level  <<Elementary physics textbook>>,
    which was used  by several generations of Russian natural scientists, engineers, and teachers.}, zoologist and ecologist
 Aleksandr Formozov, philosophers Semyon Frank and Vasily Rozanov, and
  chemist Alexey Favorsky.

\sm

 {\bf \punct  Kharkov University, 1877-1905.} After graduating from
  the gymnasium, 
 Alekseevsky entered the Natural Science Department of the Faculty of Physics and Mathematics of  Kharkov University. Two years later, he transferred
 to the Mathematical Department. In 1881, `for family reasons', he left
 the university%
 \footnote{`{\it The main obstacle to his research activity was his difficult
  financial situation'} as was written in \cite{Bah}; this is also clear from Sintsov
  \cite{Sin}). However, Sintsov in the obituary
  does not provide any details (and they are absent in other sources).}.
   In 1884, he passed examinations without attending lectures
 and graduated from the  university. At that moment he had three
 published papers \cite{Ale84-1}-\cite{Ale84-3} on the integrability of
 certain ordinary
 differential equations. They were 
  considered as a graduate (`candidate') work%
  \footnote{There is no information about his scientific adviser(s). In 
  \cite{Ale84-1}-\cite{Ale84-3}, he mentioned the names of Vasily Imshenetsky (1832-1892), who worked at the Kharkov
  University
  1877-1881 (and moved to St.-Petersburg afterward), and Aleksey Letnikov
   (1837—1888),
  who worked in Moscow.}.

 Since 1881, Alekseevsky taught mathematics, mechanics,
and physics in Kharkov in a real school%
\footnote{A  high school oriented to natural sciences, engineering,
commerce, mathematics.} (until 1891) and
 in the Second women's gymnasium (until 1887).
In 1887, he received a scholarship to prepare for a professorship
(analogous to becoming a postgraduate student in the present). 

Alekseevsky's  main work 
{\it On functions similar to the gamma function}
  was published in {\it Communications of Kharkov
Mathematical Society}%
\footnote{A well-known journal between 1879-1918. Russian: Сообщения и протоколы заседаний Математического общества при Императорском Харьковском университете; old Russian orthography: Сообщен\ия и протоколы зас\едан\ий Математическаго общества при Императорскомъ Харьковскомъ университет\е;  French: Communications et procus-verbaux de la soci\'et\'e math\'ematique de Kharkow; German: Sammlung der Mitteilungen und Protokolle der Mathematischen Gesellschaft in Charkow.}
 \cite{Ale88}--\cite{Ale89} (in fact,
it is one paper split in the middle of a sentence).
He defined the double gamma function $G(x;\alpha)$ satisfying
two functional equations
and derived an unexpectedly large collection of its properties. We discuss this impressive
work and its genesis in Sect. 3.

In 1892, this work was presented as a dissertation%
\footnote{ A `magister dissertation' in the terminology of that time. 
 In modern terminology, it apparently corresponds to a PhD.}
\cite{Ale92}. After the defence (21.02.1893), he became a privatdozent in
 Kharkov University.

In March 1894 -- November 1895 Alekseevsky visited Leipzig, G\"ottingen, and Paris.
At that time, he published a short   paper \cite{Ale94}
in German in {\it Leipz. Ber.}.
 Apparently, his work became known due to this 
publication. 



After coming back, he decided to leave the university and  teached at
a gymnasium in  town Starobelsk,
in Kharkov governorate since '{\it university \dots could do nothing to support the scientist it highly estimated}' (a quotation of a rapport of Vladimir Steklov in \cite{Sin}, who worked at  Kharkov university in that time). Alekseevsky returned  to
the university  the next year. Starting 1898, he also worked at  Kharkov Technological Institute.

Four pages of his obituary \cite{Sin}
(and a big part of \cite{StSi},  see, also, the correspondence of Steklov and 
Aleksandr Lyapunov in \cite{Ste91}) contain 
 a description of the  efforts  of several colleagues, through 1896-1904, 
  to transfer Alekseevsky to the position of extraordinary professor.
 These included
Aleksandr Lyapunov, Vladimir Steklov, Matvey Tikhomandritsky (specialist in elliptic functions), Dmitry Grave,  Konstantin Andreev (later a professor
and  dean in Moscow),  Ludwig Struve (astronomer%
\footnote{Grandson of the famous astronomer Friedrich-Vasily Struve.}),
and Ivan Osipov (chemist). The nature of  obstacles is not clear from the published literature. Steklov (letter to Lyapunov, 14.11.1903) characterized this as an 'absurd story'.


From \cite{Sin}: \it Steklov \rm [1902] \it explains, why V.~P.'s doctoral dissertation \rm [$\simeq$ habilitation] \it was delayed -- 'he could have compiled from the scientific results  a fully satisfactory dissertation'
and expresses his conviction that `if V.~P. did not take care to acquire 
formal rights, having at his disposal all the scientific materials, this circumstance can hardly diminish the scientific merits of V.~P.~Alekseevsky \dots
But we must say even more. By now, V.~P.~Alekseevsky
has prepared for printing and is already printing an extensive course of differential and integral calculus.\dots V.~P.~Ale\-k\-se\-ev\-sky's course is not a template course adapted to the usual programs of higher educational institutions.
 It is an independent academic work, where, along with pedagogical 
 requirements, it combines a rigorous and modern approach to the study of the 
 subject' \rm (the book in question was not published).

After fifth attempt,  12.01.1904 Alekseevsky got the position 
of acting ext\-ra\-or\-di\-na\-ry professor in Kharkov.
In the next year,
he moved to Tomsk for a full professorship.

\sm

It seems doubtless (see Subsect.~\ref{ss:how} and Subsect. \ref{ss:lost})
that Aleksevsky had a continuation of his work devoted to
the multiple gamma functions and the  multiple zeta functions.
However, for unknown reasons he  has not  published these investigations.  It seems also that this was one of the sources  of his career 
difficulties (see  footnote \ref{fo:SL} below). 

\sm

 {\bf\punct  Tomsk Technological Institute, 1905--1916.}
 The institute was founded in 1900 by chemist-technologist
 Efim Zubashev (see \cite{Gal}, \cite{Gal1})
  and quickly became a center of 
 education and science in Siberia. For instance, among its
 professors in the beginning of the XX century were
  Theodor Molien --- one of the founders of representation theory
  of finite groups (together with George Frobenius and William Burnside)
  and theory of associative algebras%
  \footnote{See Bourbaki \cite{Bou}, Section `Noncommutative algebra', van der Waerden \cite{Wae},
  and selected works of Molien in \cite{Mol}.} --- and the famous
   geologist Vladimir Obruchev.

According \cite{Gra2003},  Zubashev invited Alekseevsky to Tomsk for
  first time in 1900 (the Institute's the first year of work). 
In April 1905, director Zubashev offered Alekseevsky 
the position of acting ordinary professor (full professor) of theoretical mechanics 
in Tomsk.  Alekseevsky agreed and,
starting   28.11.1905,  worked in Tomsk.
He left pure mathematics, and I could not find any his scientific
publications by him
of that time.

    Granina's paper \cite{Gra2003} contains a lot of data on the educational and administrative activity of Alekseevsky in Tomsk. We mention only a few formal details.
   
    Students of  Tomsk Technological Institute were involved 
    in the tumultuous events of the First Russian Revolution
    of 1905-07.
 Zubashev (who worked under pressure from both progressists
  and con\-ser\-va\-tors)  was accused of being loyal to the striking students and
 was exiled from Tomsk. Later, he returned to the director position
 but had to leave it `for health reasons' (November 1907).  
  The institute council elected Alekseevsky  as the new director
  (rector), and
 his administrative activity was considered in local annals as
 quite successful and fruitful.

In 1911,  there was a  wave of student and
professor protests in Russia
 against the education reform  initiated by  minister Lev Kasso. The Institute was involved in these
protests. Several professors  were dismissed for oppositionism%
\footnote{The list of fired professors includes Molien and Obruchev.} (8.06.1911). Aleksevsky and the institute council tried to reinstate them. 
As a result, by an imperial edict, Alekseevsky
was ousted from the director position. However, he continued his full
 professorship and some administrative activity.

In 1916
he  stubbed his toe, didn't treat it, and died of gangrene on 13 May.

\sm




\section{ Alekseevsky's work, 1888-89}

\COUNTERS

 Alekseevsky's
work 
{\it On functions similar to the gamma function}
  was published in {\it Communications of Kharkov
Mathematical Society} %
in 1888-1889 in two papers  \cite{Ale88}, \cite{Ale89}.

\sm

{\bf\punct The  function $G(x)$.} The  Alekseevsky's paper
begins with the words: {\it The first problem, whose solution is
the topic of the present investigation, is to examine the  properties of 
a function $G(x)$ satisfying the equation:
\begin{equation}
G(x+1)=\Gamma(x)\, G(x)
\label{eq:Delta}
\end{equation}
under the condition }
\begin{equation}
G(1)=1.
\label{eq:1}
\end{equation}

 There are no references, no names, no motivations, no acknowledgements here. It is hard
to believe that the  author simply decided to attack such an exotic problem and successfully solved it. 
In reality,   the presentation of the paper was backwards to the real way of thinking; this is discussed below in Subsect. \ref{ss:how}.

Taking the logarithm of both sides  of \eqref{eq:Delta}, we come to the
difference equation
$$
\ln G(x+1)-\ln G(x)=\ln\Gamma(x).
$$
According to the Malmsten formula, see for example \cite{Bei}, formula 1.9(1),
$$
\ln \Gamma(x)= \int_0^\infty \frac{e^{-u}\,du}{u}\Bigl\{ (x-1) - \frac{1-e^{-(x-1)u}}{1-e^{-u}}\Bigr\}.
$$
Denote the integrand by $L_1(x,u)$. Alekseevsky (\S1) writes the simplest
function $L_2(x,u)$ such that
\begin{equation}
L_2(x+1,u)-L_2(x,u)=L_1(x,u)
\end{equation}
and defines $\ln G(x)$ as $\int_0^\infty L_2(x,u)\,du$,
\begin{equation}
\ln G(x):=\int_0^\infty \frac{e^{-u}\,du}{u}\Bigl\{\frac{(x-1)(x-2)}{2} - \frac{x-1}{1-e^{-u}}+ \frac{1-e^{-(x-1)u}}{(1-e^{-u})^2}\Bigr\}.
\label{eq:representation}
\end{equation}
Clearly, $G(x)$ satisfies \eqref{eq:Delta} and
\begin{equation}
G(n+1)= \frac{1^1 2^2 \dots n^n}{(n!)^n}=1^n \cdot 2^{n-1}\cdot (n-1)^1. 
\label{eq:tsel}
\end{equation}

Clearly, a solution of the functional equation \eqref{eq:1}
is determined up to a periodic factor. Further results show that 
it is a `correct' solution in an informal sense.

In \S5, using the integral representation
\eqref{eq:representation}, he gets a decomposition of 
$G(x)$ into a Weierstrass product:
\begin{equation}
G(1+x)=(2\pi)^{x/2}  e^{-\frac{x(x+1)}2-\frac\gamma 2 x^2} \prod_{n=1}^\infty \Bigl\{\Bigl(1+\frac x n\Bigr)^n e^{-x+x^2/(2n)} \Bigr\}, 
\label{eq:wei}
\end{equation}
where 
$$
\gamma:= \lim_{n\to \infty} \Bigl(1+\frac 12+\dots \frac 1n- \ln n\Bigr)
$$
is the Euler constant (cf. the similar expression for $\Gamma(x)$
in \cite{Bei}, 1.1 (3)). Next, he notices that this formula gives
an expression for $G(x)$ on the whole complex plane. 

In \S2, Alekseevsky gets 
a counterpart of the Euler definition of the gamma function, cf. \cite{Bei}, 1.1 (2),
\begin{equation}
G(x)=\lim_{n\to\infty}
\Bigl\{ (n+1)^{\frac{(x-1)(x-2)}2 } n!^{x-1}
\prod_{k=0}^{n-1}\frac{\Gamma(k+1)}{\Gamma(k+x)}
\Bigr\}.
\label{eq:euler-c}
\end{equation}

The counterpart of  the Legendre duplication formula (cf. \cite{Bei}, 1.2(15))
is (\S9):
$$
G(2x)=\frac 1{G(1/2)^2} 2^{(x-1)(2x-1)}\, \pi^{-x}\,
\Gamma(x)\, G^2(x)\,G^2(x+\tfrac 12).
$$
After complicated calculations with functions
$G(x;\alpha)$ in \S21 (see below Subsect. \ref{ss:how})
he gets the multiplication formula (cf. \cite{Bei}, 1.2(11))
in the form
\begin{equation}
G(nx)=\frac
{   n^{\frac{(nx-1)^2}2}    }
{   (2\pi)^{\frac{(n-1)(nx-1)}{2}}   }
\,\,
\prod_{j=1}^{n-1}\frac{\Gamma\bigl(x+\frac{j-1}{n}\bigr)^{n-j}}
{\Gamma\bigl(\frac jn\bigr)^{n-j}}\,\,
\prod_{j=0}^{n-1}\frac{G\bigl(x+\frac jn\bigr)^n}{G\bigl(\frac{1+j}{n}\bigr)^n}.
\label{eq:multiplication}
\end{equation}
The multiplication formula admits a simpler form obtained by Kinkelin
\cite{Kin}, see formula \eqref{eq:Kinkelin-nx}  below.

In \S11, Alekseevsky  establishes the following counterpart of the 'Stirling' asymptotic
expansion (cf. \cite{Bei}, 1.18(1)):
\begin{multline*}
\ln G(x+1)\sim
\ln \frac{\pi^{1/6} G(\frac 12)^{\frac23} }
{2^{\frac 1{36}}}+
\frac x2\ln 2\pi+ \Bigl(\frac{x^2}{2}-\frac 1{12}\Bigr)\ln x -\frac 34 x^2
+\\+
\sum_{n=1}^\infty \frac{(-1)^n B_{2n+2}}{4n(n+1)x^{2n}},
\qquad x\to+\infty.
\end{multline*}
where
$B_{2n}$ are the Bernoulli numbers, $\frac x{e^x-1}= \sum_n \frac {B_n}{n!} x^n$.

\sm

{\bf\punct Integrals and products.} Next, he evaluates
some integrals in terms of the $G$-function,
\begin{equation}
\int_0^a \ln \Gamma(x)\,dx=- \ln G(a)+(a-1)\ln \Gamma(a)-\frac{a(a-1)}{2}+
\frac a2 \ln 2\pi,
\label{eq:K1}
\end{equation}
this identity  (\S10) can be regarded as an integral representation of $\ln G(x)$.

 Keeping in mind the reflection
formula for the gamma function, we come (\S14) to
\begin{equation}
\label{eq:ln-sin}
\int_0^x \ln \sin\pi x\, dx= x\ln \frac{\sin \pi x}{2\pi}
+\ln\frac{G(1+x)}{G(1-x)}.
\end{equation}
This allows us to get some other integrals as
\begin{equation}
\label{eq:cot}
\int_0^x \pi x \cot\pi x\, dx=x\ln 2\pi+\ln \frac{G(1-x)}{G(1+x)}.
\end{equation}

We also mention the following formula (\S5),
$$
\prod_{k=0}^{n-1}\frac{G(a-e^{2\pi i k/n}x)}{G(a)}=
\prod_{m=0}^\infty \Bigl(1-\frac{x^n}{(a+m)^n}\Bigr)^{m+1}.
$$


{\bf\punct The multiple gamma functions.%
\label{ss:Gn}}
Denote 
$$G_1(x):=\Gamma(x), \qquad G_2(x):=G(x).$$
In \S16, he introduces  functions $G_n(x)$ satisfying
the relations
\begin{equation}
G_n(x+1):=G_{n-1}(x)\, G_n(x),\qquad G_n(1)=1.
\end{equation}
Namely,
\begin{equation}
\ln G_n(x)=\int_0^\infty \frac{e^{-u}du}{u} P_n(x),
\end{equation}
where
\begin{equation}
P_{n}(x+1)=\sum_{i=0}^n(-1)^i\frac{[x]_{n-i}}{(n-i)!}
\frac{1}{(1-e^{-u})^i}+(-1)^{n+1}\frac{e^{-xu}}{(1-e^{-u})^n},
\label{eq:representation2}
\end{equation}
and
$$
[x]_k:=x(x-1)\dots (x-k+1).
$$
The expressions $P_n(x)$ satisfy the difference equation
$$
P_{n+1}(x+1)-P_{n+1}(x)=P_n(x),
$$
which implies the desired property.

Alekseevsky does not develop a further theory of such functions, but
in \S17-19 he
discusses some generalities of function theory related to the Weierstrass
factorization theorem. Namely, he considers
arbitrary products of the form
$$
f_0(z)=\prod_{k=0}^\infty \Bigl(1+\frac{z}{a_k}\Bigr), 
\qquad \text{where $\Re a_k>0$, $\sum|a_k|^{-1}<\infty$}.
$$
He represents $\ln f_0(z)$ as a Laplace transform, and, applying the same trick
with functions $P_n(z)$,
constructs a sequence of functions $f_1(z)$, $f_2(z)$, \dots such that
$$f_n(z+1)= f_n(z) f_{n+1}(z).$$
He writes explicit representations of functions $f_n$ as products
of functions $G_m$ with given $m\le n$. 

\sm

{\bf \punct Double gamma functions with two periods.}
In \S20, Alekseevsky writes: {\it Functions similar to $\Gamma(x)$
considered above are special cases of more general functions. Since they are related to theta-functions, we  derive some properties of one  such function.}

Fix a parameter $\alpha $ (which initially is positive real,
 but later it becomes complex). Ale\-k\-se\-e\-vsky considers a function%
 \footnote{He uses the notation $H(x)$ or $H(x,\alpha)$.}
 $G(x;\alpha)$ such that
 $$
 G(x+1;\alpha)=\Gamma\Bigl(\frac x\alpha\Bigr)\,  G(x;\alpha),
 \qquad G(1;\alpha)=1.
 $$
 This function is defined by
 \begin{multline}
 \ln G(x;\alpha)=\\=
 \int_0^\infty \frac{e^{-\alpha u}du}{1-e^{-u}}
 \Bigl\{
 \frac{(x-1)(x-2\alpha)}{2\alpha}-\frac{x-1}{1-e^{-\alpha u}}+
 \frac{e^{-(1-\alpha) u}- e^{-(x-\alpha) u}}{(1-e^{- u})(1-e^{-\alpha u})}
 \Bigr\}
 \label{eq:representation3}
 \end{multline}
It turns out (this is non-obvious) that
$G(x;\alpha)$ satisfies the second functional equation
$$
G(x+\alpha;\alpha)=(2\pi)^{\frac{\alpha-1}{2}} \alpha^{-\frac{2x-1}{2}}
\Gamma(x)\,G(x;\alpha).
$$
Alekseevsky finds the main properties (cf. \cite{AK}) of this function, in particular,
$$
G\Bigl(x;\frac 1\alpha\Bigr)=\frac{G(\alpha x;\alpha)}{G(\alpha,\alpha)^x}
\,
\alpha^{\frac {(x-1)(\alpha x-2)}{2}},
$$
see \S21. 
In fact (\S20),
$$
G(\alpha,\alpha)=\alpha^{-\frac12} (2\pi)^{\frac{\alpha-1}{2}}
.$$
 Next (\S21),
\begin{multline*}
\frac{G(x;\alpha)}
{G\bigl(\frac{x}{\alpha+1};\frac{1}{\alpha+1}\bigr)\,
G\bigl(\frac{x}{\alpha+1};\frac{\alpha}{\alpha+1}\bigr)}
=
G(\alpha,\alpha)^{\frac{x}{\alpha+1}}
\,\alpha^{\frac{-x(x-\alpha-1)}{2(\alpha+1)}}
(\alpha+1)^{\frac{x-\alpha-1}{\alpha+1}}
\Gamma\Bigl(\frac{x}{\alpha+1}\Bigr).
\end{multline*}
If $m$, $n$ are positive integers, then $G\bigl(x,\frac mn \alpha\bigr)$
can be expressed in terms of $G(x;\alpha)$. Namely,
$$
G\Bigl(x,\frac mn \alpha\Bigr)=
\frac{n^{\frac{(x-1)(nx-m\alpha)}{2m\alpha}}}
{(2\pi)^{\frac{(n-1)(x-1)}{2}}}\,\,
\prod_{k=0}^{n-1} \prod_{j=0}^{m-1} \,
\frac{G\bigl(\frac{nx+nj+km\alpha}{mn};\alpha\bigr)}
{G\bigl(\frac{n+nj+km\alpha}{mn};\alpha\bigr)},
$$
see \S21.
This formula has numerous interesting special cases, for instance,
we can set $m=1$, or $n=1$, or $\alpha=1$, or $\alpha=1/n$, or $\alpha=n$,
or $m=n$.
The multiplication formula \eqref{eq:multiplication} also is a non-obvious
corollary of the last formula.

In \S 22, he obtains 
two limit expressions for $G(x;\alpha)$ extending \eqref{eq:euler-c},
\begin{align*}
G(x;\alpha)&=\lim_{n\to\infty}
\Bigl(1+\frac n\alpha\Bigr)^{\frac{(x-1)(x-2\alpha)}{2\alpha}}
\Gamma\Bigl(1+\frac n\alpha\Bigr)^{x-1}\,\,
\prod_{k=0}^{n-1}\frac{\Gamma\bigl(\frac{1+k}{\alpha}\bigr)}{\Gamma\bigl(\frac{x+k}{\alpha}\bigr)}\\
&=
\lim_{n\to\infty}
\alpha^{n(x-1)}
\Bigl(\frac 1\alpha+n\Bigr)^{\frac{(x-1)(x-2\alpha)}{2\alpha}}
\Gamma\Bigl(\frac 1\alpha+n\Bigr)^{x-1}\,\,
\prod_{k=0}^{n-1}\frac{\Gamma(1+k\alpha)}
{\Gamma(x+k\alpha)}.
\end{align*}
Also, he obtains
a Weierstrass type decomposition (\S22),
\begin{multline}
G(x;\alpha)=e^{ax+bx^2} \frac x\alpha
\times \\\times 
\prod_{m\ge 0, n\ge 0; (m,n)\ne (0,0)}
\Bigl(1+\frac x{m+n\alpha}\Bigr)
\,\exp\Bigl\{-\frac x{m+n\alpha}+\frac {x^2}{(m+n\alpha)^2}\Bigr\},
\label{eq:weierstrass}
\end{multline}
where the constants $a$, $b$ are defined by
\begin{align}
a:&=\Bigl(\frac {\partial}{\partial x} \ln G(x;\alpha)+\frac 1\alpha \frac{d}{dx} \ln \Gamma(x)\Bigr)\biggr|_{x=1};
\label{eq:a}
\\
b:&= \frac 12\Bigl(\frac {\partial^2}{\partial x^2} \ln G(x;\alpha)+ \frac 1{\alpha^2} \frac{d^2}{dx^2} \ln \Gamma(x)\Bigr)\biggr|_{x=1}.
\label{eq:b}
\end{align}

{\bf \punct How was the double gamma function  invented?
\label{ss:how}}
After formula \eqref{eq:weierstrass}, Alekseevsky writes:
{\it Taking this product as a definition of $G(x,\omega)$ and setting
$$
x=\frac z\omega, \qquad \alpha=\frac{\omega'}{\omega},
$$
where $\omega$, $\omega'$ are periods of elliptic functions, we get
\begin{multline}
G\Bigl(\frac z\omega; \frac{\omega'}{\omega}\Bigr)=
e^{a\frac{z}{\omega}+b\frac{z^2}{\omega^2}}
\frac{z}{\omega'}
\times \\\times
 \prod \Bigl(1+\frac{z}{m\omega+n\omega'}\Bigr)
  \exp\Bigl\{-\frac{z}{m\omega+n\omega'}+\frac{z^2}{(m\omega+n\omega')^2}\Bigr\}.
  \label{eq:pre-sigma}
\end{multline}
It is clear that this function consists of a {\rm quarter} of the factors
entering to the $\varsigma_3$-function of Weierstrass%
\footnote{The Weierstarss function, see for example \cite{Akh}, \S 13, is given by a product \eqref{eq:pre-sigma}
over $(m,n)\in \Z^2\setminus (0,0)$. In our case \eqref{eq:weierstrass},
 the product is given over
$m\ge 0$, $n\ge 0$, $(m,n)\ne(0,0$.} (or $\vartheta_1$
of Jacobi), herein it comes from $\varsigma_3$  the same way
as $\frac 1{\Gamma(z)}$ comes from $\sin\pi z$. However, in a two-periodic function,
the ratio $\frac{\omega'}\omega$ can be imaginary or complex
but for the function $G$ this condition is 
not necessary%
\footnote{i.e., a ratio of periods can be  positive real.}.

The function $G$ is analogous to the function of Heine%
\footnote{Footnote from \cite{Ale88}: {\it Handbuch der Kugelfunctionen, p. 109, or Crelle J. t 34. p. 290}, see \cite{Hei}, \cite{Hei1}. I do not understand these references (they are taken from \cite{App}). However, the attribution of this expression to Eduard Heine is correct,
see the next subsection.}, which
is composed from a {\it half} of the factors of $\vartheta_3$.
For a collation, we present an expression from Appell's memoire%
\footnote{Footnote from Alekseevsky: {\it Mathematische Annalen. B. 19. p. 84. }, see
\cite{App}},
$$
O(x,\omega,\omega')=
e^{ax^2+bx+c}\, x\prod_{m\in \Z, n\ge 0, (m,n)\ne 0}
\Bigl(1-\frac{x}{n\omega_1-m\omega_2}\Bigr) e^{\frac{x}{n\omega_1-m\omega_2}
+
\frac{x^2}{(n\omega_1-m\omega_2)^2}}
$$
for
\begin{align*}
m=0,\, \pm 1, \, \pm2,\, \dots, \pm \infty;
\\
n=0,\, +1,\,+2, \dots, +\infty,
\end{align*}
excluding the combination $m=0$, $n=0$.}

\sm

This explains the origin of this work. Looking 
to the formula in Appell's paper, Alekseevsky understood that 
it is possible to manipulate with Weierstrass-type product over a quadrant
in the lattice $\Z^2$.  However, he found
that the integral representation \eqref{eq:representation} is the key to the whole story, and the initial definition of the work became the final theorem
of his published paper (in fact, an evaluation of constants $a$, $b$ in 
\eqref{eq:weierstrass} is a highly nontrivial problem).

\sm

The work discussed here became Alekseevsky's dissertation 
\cite{Ale92}. It was published as a separate book and is a literal reprint of \cite{Ale88}-\cite{Ale89} with 
of 4 pages added to the beginning. The introduction
contains the following sentence:

\sm

{\it The known dependence between the Euler function of the second kind
{\rm[the gamma function]} and the sine {\rm[i.e., the reflection formula
$\Gamma(z)\Gamma(1-z)=\pi/\sin(\pi z)$]} made me think to extract from
the function $\vartheta_1$ of Jacobi a function similar to $1/\Gamma(x)$
and investigate its properties.

By way of analogy with the function $\Gamma(x)$, I came to  consider
 functions whose roots have the form $m\omega+n\omega'$, where $m$
 and $n$ are negative integers, and $\omega$, $\omega'$ are periods of elliptic functions.}

\medskip 

As we have seen above, Alekseevsky introduced additional pa\-ra\-me\-ters $\omega_1$, $\omega_2$ to the function
$G=G_2(x)$.
 After Barnes \cite{Bar5}, we know  that it is possible to introduce
similar parameters $\omega_1$, \dots, $\omega_n$ to the functions
$G_n(x)$ and get `{\it multigamma functions}'. Such  functions
(we use a non-standard notation) satisfy 
the collection of identities
$$
G_n(x+\omega_j; \omega_1, \dots, \omega_n)=
A_j\cdot 
G_{n-1} (x;  \omega_1, \dots,\widehat \omega_j, \dots, \omega_n)\,
G_n(x; \omega_1, \dots, \omega_n),
$$
where $A_j=A_j(x,\omega_1,\dots,\omega_n)$ are certain elementary
factors.

Actually,  Appell \cite{App} considers more general functions
$O(x;\omega,\omega_1, \dots, \omega_n)$, which are products of  certain exponentials
$$
\exp\bigl\{a_0+a_1 x+\dots +a_n x^n \bigr\}
$$
and expressions
$$
x\cdot \prod_{w} \Bigl(1-\frac x w\Bigr)
\exp\Bigl\{\sum_{j=1}^{n+1} \frac 1j \Bigl(\frac xw \Bigr)^j \Bigr\},
$$
where
$$
w= m\omega+\sum_j m_j \omega_j,
$$
and the integers $m_j$ run in the domain
$$
 m\in\Z,\qquad m_j\ge 0,\qquad (m,m_1,\dots, m_n)\ne(0,0,\dots,0). 
$$

It is strange that Appell himself did not introduce the double gamma.
In any case,  Alekseevsky's response to  Appell's paper was completely natural.

Certainly, Alekseevsky kept in mind the  multiple gamma function
with several parameters. The introductory part of 
 dissertation \cite{Ale92} contains the following 'thesis':
 
\sm 
 
 {\it 2. Roots of functions similar to $\Gamma$ have the form
 $$
 m_0+m_1 a_1+\dots + m_n a_n,
 $$
 where $m_k$ are negative integers, $a_k$ are parameters
 of the function $H$, positive integers or complex quantities.
 }

This 'thesis' had no continuation in the text of the dissertation.

\sm

Visiting Leipzig in 1894,
Alekseevsky published a short paper \cite{Ale94} in German in {\it Leipz. Ber.} (communicated
by Sophus Lie), and his work apparently became known due to this publication. The paper is a sketch of \cite{Ale88}-\cite{Ale89} with some additional claims at the end of the paper. In particular,

\sm 

{\it
Daher folgt:

\sm

(1)%
\footnote{These numbers are my addition.}
Functionen sind construirt aus einer Function $H$, die analog der
Function $\Gamma$ ist.

\sm

(2) Dieser Satz giebt auch eine \"ahnliche Beziehung der
Modulfunctionen von Hermite zur Function $H$.

\sm

(3) Aus der Betrachtung der Functionen h\"oherer Ordnung, die
analog der Function $\Gamma(x)$ sind, folgt:

\sm

 Die Euler'schen Functionen von Herrn Appell sind Producte
von Functionen, die analog der Gammafunction sind.}

\sm

Clearly, in (3) Alekseevsky claims that an Appell function
is a product of two multiple gammas%
\footnote{In (1) 'Heine function' is \eqref{eq:O}; it is the $q$-gamma
function up a normalization,
which can be expressed in $G(x;a)$ by the reflection formula
(see my next subsection).
Since we can express theta as a product of 4 double gammas, we
can express any double periodic function, see for example \cite{Akh}, \S 19.
\newline
But I do not understand the meaning of '{\it Modulfunctionen von Hermite}'
and the claim (2).}.

\sm

In any case, the multiple gamma function with several parameters had not
appear in  Alekseevsky's published papers  and actually was introduced
by Barnes \cite{Bar5}.
 In Subsect. \ref{ss:lost}, we continue this discussion.

\sm

{\bf \punct  The reflection formula.}  In particular, the Leipzig paper
\cite{Ale94} contains a counterpart
of the reflection formula for the gamma function,
\begin{equation*}
G(1+x;\alpha)\, G(-x;-\alpha)= C\cdot O(x),
\end{equation*}
where 
\begin{equation}
O(x):=\prod_{k=1}^\infty (1-q^{2k} e^{2\pi i x}), \qquad
\text{where $q=\pi i\alpha$,}
\label{eq:O}
\end{equation}
(this follows from  the previous subsection's reasonings).
Applying the Jacobi triple identity (see, for example \cite{Akh}, \S 18), we can easily express
the standard theta functions in the double gamma.

\sm










\sm

{\bf \punct A counterpart of the $\psi$-function.} Let us return to
the work \cite{Ale88}-\cite{Ale89}.
 Recall that
the classical $\psi$-function (see, for example \cite{Bei}, 1.7) is defined by
$$
\psi(z)=\frac{d}{dz}\ln\Gamma(z)=\frac{\Gamma'(z)}{\Gamma(z)}=
\int_0^\infty \Bigl(\frac{e^{-t}}{t}-\frac{e^{-zt}}{1-e^{-t}}\Bigr)\,
dz.
$$
In \S3 Alekseevsky introduces a similar function
$$
\phi(x)=\frac{d}{dx}\ln G(x)
$$
(in fact, both functions appeared above in \eqref{eq:a}, \eqref{eq:b}).
The function $\phi$ can be expressed explicitly in terms of $\psi$,
$$
\phi(x)=(x-1)(\psi(x)-1)+\phi(1), \qquad \phi(1)=-\tfrac12+\tfrac12 \ln 2\pi.
$$
for $a=1$ this gives \eqref{eq:wei}. So, here he has not come to a new special function, but  he has found some interesting formulas.

Expanding $e^{ux}$ in \eqref{eq:representation} in $(ux)^n$,
we can get (\S7) an expansion of $\ln G(x+a)$ in a power series,
\begin{equation*}
\ln G(x+a)=\ln G(a)+x\phi(a)+\frac{x^2}2\phi'(a)+
+\sum_{j\ge 3}(-1)^{j+1}\frac{x^j}{j} C_j,
\end{equation*}
where 
$$
C_j=\sum_{k=1}^\infty \frac{k}{(a+k-1)^j}.
$$
For $a=1$ we have $C_j=\zeta(j-1)$, for general $a$ we have an expression
in terms of Hurwitz zeta $\zeta(a;j)$ and $\zeta(a;j-1)$.

Alekseevsky obtains some other interesting identities for the function
$\phi$, for instance,
$$
\phi(1+x)=-\tfrac 12+\tfrac 12\ln 2\pi-x(1+\gamma)+
\sum_{k=1}^\infty\frac{x^2}{k(x+k)},
$$
and
\begin{equation*}
\phi(a+x)=\phi(a)+x\psi(a)+
\sum_{j=1}^\infty (-1)^{j-1}
 \frac{x(x-1)\dots(x-j+1)}{j(j+1)(a)_j}.
\end{equation*}

\section{On the history of the double gamma function} 

\COUNTERS

{\bf \punct Kinkelin's work, 1860.}
Hermann Kinkelin \cite{Kin} introduced a function%
\footnote{His original notation is ${\bf G}(x)$.}
$K(x)$ satisfying the identity
\begin{equation}
K(x+1)=x^x K(x), \qquad K(1)=1
\label{eq:Kink}
\end{equation}
by (\S6)
$$
\ln K(x)= \int_0^x \ln \Gamma(x)\,dx+\frac{x(x-1)}{2}+
\frac 12 x\ln 2\pi.
$$
Comparing this expression with
\eqref{eq:K1}, we see that the double  gamma function $G(x)$ 
is equivalent to $K(x)$, namely 
$$
G(x)\, K(x)=\Gamma(x)^{x-1}
$$ 
 as it was observed by Beaupain \cite{Bea1903}.
 
 Kinkelin obtained some properties of the function $K(x)$.
  In particular, he wrote
 its representations as limits of products and expressions of 
 elementary integrals in terms of $K(x)$ (\S\S10-12).
  The most remarkable
of his identities for $K(x)$ is the multiplication formula, which was established in a better
form than  Alekseevsky's \eqref{eq:multiplication}:
\begin{equation}
K(nx)=
n^{\frac 12 nx(nx-1)+\frac 1{12}}
\,
 \wt \omega^{-\frac12(n^2-1)}\cdot
\prod_{j=0}^{n-1} K\bigl(x+\tfrac jn\bigr)^n, 
\label{eq:Kinkelin-nx}
\end{equation}
where $\wt \omega $ is a certain constant.

He found several expressions for this constant.
For each $n>1$,
\begin{equation}
\ln \wt \omega=\frac{n}{n^2-1}\Bigl(\frac 1{12} \frac{\ln n}n+\sum_{j=0}^{n-1}
\ln K\Bigl(\frac jn\Bigr) \Bigr),
\label{eq:prelimit}
\end{equation}
see \S8. In particular, for $n=2$ this means
$$
\ln K\bigl(\tfrac12\bigr)=\tfrac 34 \ln\wt \omega-\tfrac 1{24}\ln 2.
$$ 
So, the constant $\wt \omega$ is connected  to 
$G\bigl(\frac 12\bigr)$, which appears in many formulas for
the double gamma, by a simple relation.
Passing to the limit as $n\to \infty$ in \eqref{eq:prelimit}, we get
$$
\ln\wt \omega:= 2\int_0^1\ln K(x)\,dx. 
$$
He also  expressed $\wt\omega$
in terms of the  zeta-function $\zeta(x)=\sum_{n=1}^\infty n^{-x}$,
$$
\tfrac12 \ln \wt \omega=-\tfrac 1{24}+\frac 13 \gamma
+\sum_{\lambda=1}^\infty
 \frac{\zeta(2\lambda+1)-1}{(2\lambda+1)(2\lambda+3)},
$$ 
where $\gamma$ is the Euler constant (see \S8).
Another expression is given by
 the following counterpart (\S9)  of the Raabe 
formula (see \cite{Bei}, 1.9(19)) for the gamma function:
$$
\int_x^{x+1} K(t)\,dt=\tfrac 12 \ln \wt \omega+\tfrac 14 x^2(2\ln x-1)
$$
(Alekseevsky later repeated this formula).

Kinkelin also establishes the asymptotics (\S8):
$$
K(n+1)=\prod_{j=1}^n j^j \sim \wt w^{1/2} \,
e^{-\frac 14 n^2+\frac 1{12}} n^{\frac 12(n^2+n)+\frac 1{12}},
\qquad n\to+\infty.
$$

The expression $A:=(\wt \omega)^{1/2}e^{1/12}$  is now called the {\it Glaisher--Kinkelin constant}, cf. Glaisher \cite{Gla}.

\sm

 Kinkelin's work didn't attract immediate attention and became a subject of interest only in the second half of the 1890s.
 The problem \eqref{eq:Kink} seems even more artificial than
 Alekseevsky's \eqref{eq:Delta},
 so it is natural to ask
about his way to the 'Kinkelin function'.
Apparently, it was a result of the following 'speleology'.

 \sm

For a function $f$, denote 
$$
S_n(x):=\sum_{j=0}^{n-1} f\Bigl(x+\frac jn\Bigr).
$$
In his previous work, 
Kinkelin \cite{Kin0} tries to find all functions
$f(x)$ satisfying the 'addition theorem' in the following sense:
for each $n>0$ there are constants $\phi(n)$, $\psi(n)$ such that
\begin{equation}
S_n(x)=\psi(n)\cdot n  f(nx)+\phi(n).
\label{eq:addition}
\end{equation}

Apparently, this work was a response to the paper by 
Joseph Raabe \cite{Raa} on Bernoulli polynomials $B_p(x)$.
Recall that they are defined by the generating function:
\begin{equation}
\sum_{p=0}^\infty B_p(x)\frac{t^p}{p!}= \frac{t e^{xt}}{e^t-1}.
\label{eq:bern-pol}
\end{equation}
Their main property is the difference equation
$$
B_p(x+1)-B_p(x)=p x^{p-1}.
$$
Raabe showed that they satisfy the addition theorem of the type
\eqref{eq:addition}: 
$$
\sum_{j=0}^{n-1} B_p\Bigl(x+\frac j n\Bigr) =
n^{1-p} B_p(nx) 
$$
(this easily follows from \eqref{eq:bern-pol}).
In fact, Raabe used another normalization of the Bernoulli
polynomials, $\wt B_p(x)=B_p(x)-B_p(0)$,
 and his $\phi(n)$ differed from 0).

Kinkelin examines other solutions of the  problem \eqref{eq:addition}  
(I am not sure that his proof is precise). There are also  solutions
$$
f_q(x)=\frac {(-1)^q}{(q+1)!}\,\frac{d^{q+1}}{dx^{q+1}}\ln \Gamma(x).
$$
Notice that now $f_q(x)$ satisfies  the difference equations
$f_q(x+1)-f_q(x)=\frac 1{x^{q+1}}$. This family also was not `new'
in Kinkelin's time.

However, he tries to vary the condition
\eqref{eq:addition}. The function $f(x)=\ln \Gamma(x)$
satisfies the addition theorem
\begin{equation}
\sum_{j=0}^n f\bigl(x+\tfrac jn\bigr)=
f(nx)+ \frac{n-1}2\ln(2\pi)+\bigl(\tfrac12-nx \bigr)
 \ln n,
 \label{eq:k1}
\end{equation}
which is equivalent to the multiplication formula for the gamma function.
Also, it satisfies  the difference equation
\begin{equation}
f(x+1)-f(x)=\ln x.
\label{eq:k2}
\end{equation}
We take primitives from both sides of identities \eqref{eq:k1},
\eqref{eq:k2} and get two equations
 for the primitive  $F(x):=\int f(x)\,dx$.  Equation \eqref{eq:k2}
 yields
$$
F(x+1)-F(x)=x\ln x-x+C.
$$
Correcting $F(x)\mapsto \wt F(x)$ by an appropriate quadratic summand, we can reduce the last equation to the form
$$
\wt F(x+1)-\wt F(x)=x\ln x,
$$
i.e., the equation \eqref{eq:Kink} for $\ln K(x)$.
The integration of \eqref{eq:k1} gives an addition theorem
for $\wt F(x)$ or a multiplication theorem
for the Kinkelin function $K(x)$.

\sm

{\bf \punct The double gamma function, 1899-1907.}
The work of Alekseevsky became known (perhaps due to the publication
\cite{Ale94}). After this, the double gamma function was
 a topic of publications
 by Ernest Barnes \cite{Bar1}, \cite{Bar2}, \cite{Bar3}, \cite{Bar4},
 \cite{Bar-cap}, \cite{Bar5}, Jean Beaupain \cite{Bea1903}, \cite{Bea1904}, \cite{Bea1907},
Godfrey  Hardy \cite{Har1}, \cite{Har2}, and
 Vladimir Steklov \cite{Ste}. The most important result of these works
  was Barnes's introduction
 of the multiple gamma functions \cite{Bar5}.

  In particular, the third Barnes's
publication  \cite{Bar1} and his further works
\cite{Bar2}--\cite{Bar-cap} were devoted to the double gamma.
In \cite{Bar-cap} he found the first application: 
 he expressed the  electrostatic capacity coefficients 
of a pair 
of non-intersecting spheres in  terms of the double gamma.

In the work \cite{Bar5} he introduced a  new topic, namely the {\it multiple
gamma function}; it is a function $G_n(x)$ depending on parameters 
$\omega_1$,
\dots, $\omega_n$ and satisfying functional equations
$$
G_n(x+\omega_j; \omega_1, \dots, \omega_n)=
A_j\cdot 
G_{n-1} (x;  \omega_1, \dots,\widehat \omega_j, \dots, \omega_n)\,
G_n(x; \omega_1, \dots, \omega_n),
$$
where $A_n$ are certain elementary expressions and $G_{n-1}$ is the multiple gamma of the previous level.

Barnes became famous due to his later works
of 1907--08 on the 'Mellin--Barnes integrals', i.e.,
integrals  of the  type
$$
\int_{-i\infty}^\infty z^{-s} 
\frac{\prod \Gamma(a_k+s)\prod\Gamma(b_l-s)}
{\prod\Gamma(c_m+s)\prod\Gamma(d_n-s)}\,ds.
$$
Integrals of such kind appeared earlier in works by Salvatore
Pincherle and Hjalmar Mellin. In any case, Barnes's research
was a breakthrough --- thanks to him such integrals
became a powerful and flexible tool of theory of special functions,
see for example \cite{Whi}, \cite{AAR}, \cite{Mar} (in particular, Mellin--Barnes integrals were the main tool of creations of tables of integrals  \cite{PBM}, see \cite{PBM-plus}; see also the impressive table of Mellin-Barnes integrals
in volume 3 of \cite{PBM}, Chapter 8). 
 
 \sm
 
{\bf \punct The higher Kinkelin functions.}
Beaupain \cite{Bea1903} introduced higher Kinkelin functions,
i.e., functions, which satisfy the identity
$$
K_n(x+1)=x^{x^n} K_n(x), \qquad K_n(1)=1.
$$ 
 Alekseevsky \cite{Ale04} showed that systems of functions
 $G_n(x)$ (see above Subsect. \ref{ss:Gn}) and $K_n(x)$ can be expressed through each other,
 \begin{equation}
\ln K_n(x)=\sum_{j=0}^n (-1)^{j} (\Delta^j x^n)\, \ln G_j(x+j), 
 \end{equation}
where $\Delta f(x):=f(x+1)-f(x)$. 
 
 \sm 
 
 {\bf \punct The gap, 1908-1975.}
 After 1907, the double gamma was forgotten.
 
As  was observed in \cite{Ada},   this function was
 mentioned in the second edition \cite{WhWa} 
 of E.~T.~Whittaker,  G.~N.~Watson  {\it A course of Modern Analysis} 
 (Chapter XII, Examples 48-50, namely \eqref{eq:wei} as a definition,
 \eqref{eq:Delta}, \eqref{eq:tsel}, \eqref{eq:cot}).
 Next, it was mentioned in the fourth edition \cite{Gra1963} of
 Gradshteyn,  Ryzhik {\it Tables of integrals, sums, and series,}
lines 6.441, 8.333, with references to  Whittaker--Watson.
 
 The double gamma was absent in other handbooks on special functions as {\it Higher transcendental functions} \cite{Bei}, or  Prudnikov,
  Brychkov, Marichev, \cite{PBM}.

I have seen (mainly, due to \cite{Gus}) a few other appearances of such functions in literature. In 1933,
L.~Bendersky \cite{Ben} published a large detailed investigation
on higher Kinkelin functions (he probably didn't know about his predecessors).

Boris Rimsky-Korsakov published a paper \cite{Rim1} on the reciprocity theorem
for higher Kinkelin functions:
\begin{equation*}
\frac
{\Bigl\{\prod\limits_{k=0}^{q-1}K_n\bigl(pz+\frac{kp}{q}\bigr)\Bigr\}^{q^n} }
{\Bigl\{\prod\limits_{k=0}^{p-1}K_n\bigl(qz+\frac{kq}{p}\bigr)\Bigr\}^{p^n}}=
\Bigl(\frac p q\Bigr)^{\frac{p^n}{n+1}\sum_{k=0}^{n-1} B_{n+1}
\bigl(qz+\frac{kq}p\bigr)}
\pi_n^{(q^{n+1}-p^{n+1})/2},
\end{equation*}
where $B_m(\cdot)$ are the Bernulli polynomials and $\pi_n$ are certain 
constants. For $n=1$ this statement is 
 the Schobloch reciprocity formula (\cite{Sch}, \cite{Sri})
  for the gamma function.
He also published  the paper \cite{Rim2} and PhD thesis%
\footnote{I have not seen it.} \cite{Rim3} about further generalizations of the gamma function. 

The double gamma function was discussed in the obituary of Barnes, written
by Edmund Whittaker \cite{Whi}. 

\sm

{\bf \punct A new appearance.} The double gamma function appeared again
in the work of Takuro Shintani \cite{Shi1}-\cite{Shi2}.
After that, the function started showing up in various math problems,
I only present a few references \cite{Vig}, \cite{Sar}, \cite{Var},
\cite{NY}, \cite{Ada}, \cite{Ost}, \cite{AK}, \cite{KK}.

Shintani \cite{Shi2} also introduced  the {\it double sine}
\begin{multline*}
\cS_2(x;\omega_1,\omega_2)
=(2\pi)^{\frac{\omega_1+\omega_2-2x}{2\omega_1}}
\omega_2^{\frac{x(\omega_1+\omega_2)}{\omega_1\omega_2}}
\frac{G\bigl(\frac x{\omega_1};\frac{\omega_2}{\omega_1}\bigr)}
{G\bigl(1-x+\frac{\omega_2}{\omega_1};\frac{\omega_2}{\omega_1} \bigr)}
=\\=
\exp\Bigl\{
\int_0^\infty\Bigl(
\frac{\sinh\bigl(z-\frac12(\omega_1+\omega_2)\bigr)}
{\sinh\bigl(\frac12\omega_1 t\bigr)\,\sinh\bigl(\frac12\omega_1 t\bigr)}
-\frac 1{\omega_1 \omega_2 t}(2z-\omega_1-\omega_2)
\Bigr)\frac{dt} t
\Bigr\}
\end{multline*}
(a special case of the double sine
corresponding to equal periods  was considered by Otto H\"older \cite{Hol}). 
This function is  important by itself. Much later (and independently) 
Ludvig Faddeev \cite{Fad} introduced the {\it quantum dilogarithm}, which is very
similar to the double sine. However, it appeared in a completely new context,
and  the term 'quantum dilogarithm' is now even more common than 'double sine'.
 
 For instance, see some nice applications of the double sine in
rep\-re\-sen\-ta\-ti\-on-the\-o\-re\-ti\-cal spirit (which were the reason
of my interest in this topic) in  \cite{Kash}, \cite{Fock}, \cite{KLS}).

\sm

{\bf\punct Why were the double gamma and the multiple  gamma 
  forgotten 70 years?}
I'll try to speculate on possible causes. One of the reasons is clear:
the function was  different significantly from special functions that were
considered at that time.
In my opinion, there were also some accidental reasons.

Firstly, both Alekseevsky 
 and Barnes (1874-1953) left mathematics. Alekseevsky
 became a professor of mechanics and  a rector of a technological institute (his last research publication was in 1904).
 Barnes   became an Anglicanic  deacon,
a theologist, and later a  bishop of Birmingham,
two of his last research publications are dated  1908 and 1910,  see
\cite{Whi}).

Secondly, this can be partially related to the structure
of  \cite{Ale88}-\cite{Ale89}. The work can seem like a free fantasy,
since its real external
relations become clear only at the end of the second paper.
There was a chance for a new appearance of the double gamma
after a historical-mathematical paper by Viktor Gussov \cite{Gus},
which contained a good exposition
of \cite{Ale88}-\cite{Ale89} (I think  many Russian mathematicians have
seen this article and 
were surprised  reading it). But that didn't happen.

Thirdly, a meticulous reader could decide that the actual topic of 
the works on the double gamma are
the functions $\ln G(x)$, $\ln G_n(x)$, $\ln G(x;a)$. For instance,
in  Alekseevsky's papers, the double gamma itself (without  a logarithm) never appears in the right hand sides of  integrals, which contrasts with the behavior of $\Gamma(x)$. The classical gamma function appears in the right hand sides of
a huge number of  definite integrals, see for example \cite{Gra1963}, \cite{PBM}; on the other
hand,  appearances of $\ln \Gamma(\cdot)$ in the
right hand sides of integrals are
 relatively rare.
Next, integrals like \eqref{eq:ln-sin}, \eqref{eq:cot} are not too impressive
(the right hand side is not better than the left hand side)
and admit other transformations.
 Only in \S22 of \cite{Ale88}--\cite{Ale89}  a reader  can understand
  that products of the double gammas can participate in Heine's theory of $q$-hypergeometric functions%
 \footnote{Also, $q$-hypergeometric functions were almost unknown in Russia until
the
appearance of interest in representations of classical groups over
finite and $p$-adic fields in the 1960-s. So, to a Russian reader,
motivation of \cite{Ale88}-\cite{Ale89} was unclear.},

\sm

Now there are no reasons for such  doubts. 
 Ludvig Faddeev \cite{Fad}
observed that the double sine 
is a kernel of a natural integral operator in $L^2(\R)$ satisfying highly
nontrivial identities, see \cite{Fock0}, \cite{Kash}, \cite{Fock}. See also
the work \cite{KLS}, where  certain spectral measures for quantum groups
 are expressed as  products
 of the double sines. Recall that products of the usual gamma functions 
 appear as spectral densities for differential and difference operators
 (as far as they admit explicit evaluations),  such products also appear
 as spectral densities  in numerous problems
 of harmonic analysis for semisimple Lie groups.

 See, also a recent work \cite{KK}.
 
 \sm

The
 term '{\it Alexeiewsky function}' was used in several papers
 published in 1899-1907.
The widely used modern term '{\it Barnes double gamma}' was originated
by Shintani \cite{Shi1}--\cite{Shi2}, who found old Barnes works. Alekseevsky, who was mainly
the author of one big work \cite{Ale88}, \cite{Ale89}, \cite{Ale94}, \cite{Ale04},  was forgotten at that time, while
Barnes was well known  due to 
the Barnes integrals (see above).


\section{Alekseevsky: his historiography,  lost work, and publications}

\COUNTERS



{\bf\punct Historial sources and historiography.}
I used the following contemporary historical sources on the biography
of Alekseevsky:

\sm
 
 --- Obituary written by the
 president (1906-1946) of the Kharkov Mathematical Society Dmitry Sintsov \cite{Sin}.
 
\sm 

---  Steklov and Sintsov's review   \cite{StSi} of  Alekseevsky's works 
 (apparently written by Steklov) supporting Alekseevsky's application  to the position
 of extra-ordinary professor in  Kharkov University.

\sm

--- Steklov's correspondence   published in \cite{Ste91}.

\sm

As told above, Alekseevsky's work was forgotten. But his name 
several times was mentioned in  histories of mathematics
and education in Russia. 

\sm
 
 ---  Viktor Gussov  wrote an interesting and thorough
  dissertation \cite{Gus0}
    on investigations of the gamma function and the Bessel functions
in Russia  (the work was published in \cite{Gus}--\cite{Gus1}). He found Alekseevsky's works and  exposed them in \cite{Gus0}, \cite{Gus}. 
Aleksey Markushevich was an opponent of Gussov's thesis and
used it in his 'Essays' on the history of complex analysis  \cite{Mark}. However,
he did not mention Alekseevsky  in his 'Essays', nor in
his historical investigation \cite{Mark1}, nor in his famous  treatise 
{\it Theory of  functions of complex variable}. Clearly,
Gussov's paper \cite{Gus}
was used in  the publications mentioned below. 
 
 \sm 
 
 --- Some information is contained in  \`Esfir Bakhmutskaya's work
 \cite{Bah} on mathematics in  Kharkov University.
 
\sm 
 
 --- The Tomsk period of Alekseevsky's biography was well-investigated 
 in several works about
  the history of science in Tomsk (see Nikolay Krulikovsky \cite{Kru},
  \cite{Kru1},
  Yulia Granina \cite{Gra2003}, Belomestnyh \cite{Bel}).

 \sm

 \sm
 
 {\bf \punct The multiple gamma function and a lost work.%
\label{ss:lost}}
Let us return to the discussion in Subsect. \ref{ss:how}.
In their review of Alekseevsky works \cite{StSi}, Steklov and Sintsov
wrote:

\sm 

{\it  Now as well, Alekseevsky continues to successfully develop 
his theory of functions, which  he now calls gammamorphic. He came to functions
of an even more general type,  for which Kinkelin functions {\rm [Steklov keeps in mind $G_n(x)$]} are limit cases,
and which are related to the known Riemann $\zeta(s)$ and many other 
functions important in Analysis and its applications.

These recent studies were reported at the meetings of the Kharkov Mathematical Society in 1902 and 1903.  
Part of them are currently beeing printing in
our society[s Communications.
}

\sm

We see in \cite{Khar} the following titles of Alekseevsky's talks at
the Kharkov Mathematical Society:

\sm

15.02.1902. On the Riemann $\zeta$-function of three arguments.

03.12.1902. On the Riemann $\zeta$-function of two arguments.

24.10.1903. Note on the Kinkelin functions.

31.10.1903. Towards the theory of gammamorphic functions.

\sm

Clearly, the third talk corresponds to \cite{Ale04}.

 Sintsov \cite{Sin} mentions 
that an  article `{\it On Riemann function $\zeta$}'
disappeared
and writes that there is  hope to find it in   Alekseevsky's
posthumous papers.

\sm 

It seems clear that the topic of these works  
were the multiple gamma functions and
related $\zeta$-functions of the type 
$$\sum_{m,n\ge 0} (z+n \omega_1+m \omega_2)^{-s}$$
(something similar to Barnes \cite{Bar5}).

In any case, Alekseevsky  has not published anything on this topic
except the  claims in \cite{Ale89}, \cite{Ale92}, \cite{Ale94} quoted above%
\footnote{\label{fo:SL} Letter of Steklov to Lyapunov, 29.12.1902: \it
Alekseevsky, as usual, drags on and does not submit his article 
{\rm [Yu.~N.: clearly, this sentence refers to \cite{Ale04}]}, although he says that he will submit 
it `one of these days', but this `one of these days' has been going on since September.}.

\sm

{\bf\punct  Alekseevsky's publications.%
\label{ss:publications}} The list of his research publications is
 con\-ta\-i\-ned in our bibliography,
 \cite{Ale84-1}--\cite{Ale04} (all of them belong 
  to the Kharkov period of his work). 
 
\sm 
 
  Russian library catalogues
 show several editions of his lecture courses and problem books
 (in  Kharkov Technological Institute  and in  Tomsk Technological Institute):

\sm

Алексеевский В. П.  {\it
Дифференциальное исчисление : Лекции, чит. в весен. семестре 1900 г.}   Харьков: типо-лит. С.~Иванченко, 1900., 317 с. 

\sm

{\it Задачи по дифференциальному и интегральному исчислениям, 
пред\-ла\-га\-в\-ши\-е\-ся проф.~Алексеевским в Х\rm[\it арьковск.\rm] \it Т\rm[\itехнол.\rm] 
\it И\rm[\it н-те\rm ]}. - Харьков : изд. студ.-технологов Н.~Александровича и Н.~Праведникова, 1907. - [2], 59 с.

\sm

{\it Динамика сост. по лекциям и под ред. В. П. Алексеевского студ. И. Никишовым, 1912-1913 уч. г.} ; Том. технологич. ин-т императора Николая II.
Типо-лит. Сиб. Т-ва Печатн. Дела 1913.

\sm

{\it Сборник задач по динамике с подробными решениями, предлагавшихся про\-фес\-со\-ром В.~П.~Алексеевским в Томском технологическом институте в 1913-14 уч. году на III и IV сем.}  Сост. Г. Берштейн, студент Мех. отд-ния. - Томск: типо-лит. Сиб. т-ва печ. дела, 1914 , 47 p. 

\sm

Алексеевский В.П. {\it Динамика : Лекции 1915-1916 уч. г.}  Том. технол. ин-т. Томск: Т-во "Печатня С.~П.~Яковлева", 1916, 237 с.

\sm

Also, \cite{Sin} and \cite{Gra2003} contains information about some additional li\-tho\-gra\-phic editions of his lecture courses in  Kharkov University and
in  Tomsk Technological Institute.

\section{Conclusion}

The double gamma  was introduced in the work of 
Vladimir Alekseevsky in 1888-1889, he also derived
main identities involving this function
(cf. the list of identities \cite{AK}). 
His publication \cite{Ale94} had a big influence on
the early works by
Ernest Barnes, who also  introduced
and investigated the multiple
gamma function \cite{Bar5}.  The double gamma function
 was Alekseevsky's main work
and he was forgotten for the same period  that the double gamma itself was forgotten. 

In this paper, we observed 
 an interesting case: that the results of a work, which were very difficult to replicate, were forgotten for 70 years.
This happened despite the fact that some of the authors writing about 
the double gamma, were well-known and that the results were vivid.
The cause was probably an unfortunate confluence of several circumstances.

Modern mathematicians use complicated languages, 
having trouble understanding each other. Large losses of 
nontrivial information  over the intervals of several decades
are unavoidable (and is, in fact, already occuring).
 In our case, however, the papers on special
functions and complex analysis of 1880-1910s were written
in a simple language one accessible  for mathematicians of
1970-s and since.
They could have been easily read
and used.


\tt
University of Graz;

High School of Modern Mathematics MIPT; 

Moscow State University, MechMath Faculty;

Univesity of Vienna, Faculty of Mathematics.

URL: https://www.mat.univie.ac.at/$\sim$neretin/
	

\begin{thebibliography}{cc}

\bibitem[Ada2014]{Ada}
Adamchik V.S.
{\it Contributions to the Theory of the Barnes Function},
Int. J. Math. Comput. Sci. 9 (2014), no. 1, 11-30.

\bibitem[Akh1970]{Akh}
Akhiezer N. I.
{\it Elements of the theory of elliptic functions.}
American Mathematical Society, Providence, RI, 1990.



\bibitem[Ale1884.1]{Ale84-1}
	Алексеевский В. П.  {\it Об интегрировании уравнения 
	$\frac{d^ny}{dz^n}+\frac \alpha z \frac{d^{n-1}y}{dz^{n-1}}+ \beta{y}=0,
	$
} (Russian)
	   Сообщ. и протоколы заседаний Матем. общ. при Император. Харьков. унив. 1884 года,  1884, 1,  41-64;
[Alexeiewsky V. P. {\it On the integration of the equation 
	$\frac{d^ny}{dz^n}+\frac \alpha z \frac{d^{n-1}y}{dz^{n-1}}+ \beta{y}=0,
	$} Charkow Ges.,  1884, 1,  41-64]
	JFM%
	 \footnote{JFM database; the electronic version of
 {\it Jahrbuch \"uber die Fortschritte der Mathematik}, 1868-1942.
 It is contained in the {\it Zentralblatt f\"ur Mathematik
	 database.}} 16.0290.01. 

\bibitem[Ale1884.2]{Ale84-2}
	 Алексeевский В. П. {\it Заметка об обобщении уравнения Риккати.} (Russian)  Сообщ. и протоколы заседаний Матем. общ. при Император. Харьков. унив. 1884 года, 1884, 1,  80-82 [Alexeiewsky V. P.  {\it Note on the generalization of the Riccati equation.} [Charkow Ges.], 1884, 1,  80-82]  
	 JFM 16.0290.01]

  
\bibitem[Ale1884.3]{Ale84-3}	
	 Алексеевский В. П.  {\it Об интегрировании уравнения 
	 $x^2y'''+Axy'''+By'+Cx^\mu y=0,$} (Russian)
	 	  Сообщ. и протоколы заседаний Матем. общ. при Император. Харьков. унив. 1883 года, 1884, 2,  115-126
	 (Russian) [Alexeiewsky W. {\it On the integration of the equation  $x^2y'''+Axy'''+By'+Cx^\mu y=0,$}
	 Charkow Ges., 1884, 2,  115-126]
  JFM 16.0288.01. 
	
\bibitem[Ale1885]{Ale85}
	 Алексеевский В. П., {\it Об интегрировании одного линейного диффенциального уравнения $n$-го порядка} (Russian)
	 Сообщ. и протоколы заседаний Матем. общ. при Император. Харьков. унив. 1884 года, 1885, 3,  222-232. [Alexeiewsky V. P. {\it On integration of	 one linear differential equation of $n$-th order},  Charkow Ges., 1885, 3,  222-232] JFM 17.0328.01	

	
\bibitem[Ale1888]{Ale88}	
 Алексеевский В. П.  {\it О функциях подобных функции гамма} (Russian),
  Сообщ. Харьков. матем. общ. Вторая сер., 1:1 (1888),  169-192
 [Alexeiewsky W. P. (Alexejewsky W.) {\it On the functions that are similar to the gamma function.}], Charkow Ges., 1:1 (1888),  169-192]
  JFM 22.0439.03. 

\bibitem[Ale1889]{Ale89}	
 Алексеевский В. П.  {\it О функциях подобных функции гамма} (Russian),
  Сообщ. Харьков. матем. общ. Вторая сер.,
1:2 (1889),  193-238.
 [Alexeiewsky W. P. (Alexejewsky W.) {\it On the functions that are similar to the gamma function.}, Charkow Ges., 
1:2 (1889),  193-238]  JFM 25.0751.02 	
	
\bibitem[Ale1892]{Ale92}  Алексеевский В. П.
{\it О функциях подобных функции Гамма.} (Russian) Диссертация. Харьков: тип. Зильберберга, 1892,
 70 с. [Alexeevsky V.P. {\it On the functions that are similar to the gamma function.} Dissertation. Printing house of Zilberberg, Kharkov, 1892, 70p.] 	
	
\bibitem[Ale1894]{Ale94}	
	 Alexejewsky W. P. {\it Ueber eine Classe von Funktionen, die der Gammafunktion
analog sind} (German), Berichte \"uber die Verhandlungen der K\"oniglich-S\"achsischen Gesellschaft der Wissenschaften zu Leipzig, 46 (1894) 268-275.
 JFM 25.0751.02

\bibitem[Ale1895]{Ale95}
 Алексеевский В. П.  {\it Об автоморфной функции аналогичной экспонентной} (Russian),   Сообщ. Харьков. матем. общ. Вторая сер., 4 (1895),  253-262
[Alexeiewsky V. P. {\it On an automorphic function that is analogous to the exponential},
 Charkow Ges., 4 (1895),  253-262] JFM 26.0466.02


 \bibitem[Ale1899.1]{Ale99-1} 
 Алексеевский В. П.  {\it О законе взаимности простых чисел}
(Russian), Сообщ. Харьков. матем. общ. Вторая сер., 6 (1899),  200-202 
 [Alexejevsky~W.~P. {\it On the reciprocity law of primes}, Charkow Ges., 6 (1899),  200-202] JFM 30.0184.03
  
\bibitem[Ale1899.2]{Ale99-2} 
Алексеевский В. П.  {\it Об определении длины в неэвклидовой геометрии} (Russian) Сообщ. Харьков. матем. общ. Вторая сер., 6 (1899),  139-153 [Alexejewsky W. P. {\it On a definition of a length in the non-Euclidean geometry}, Charkow Ges., 6 (1899),  139-153] JFM 29.0411.03



 \bibitem[Ale1904]{Ale04}
 Алексеевский В. П. {\it Зависимость между Кинкелиновыми и гаммаморфными функциями.} (Russian) Сообщ. Харьков. матем. общ. Вторая сер., 8 (1904),  123-135. [{\it Alexejewsky W. P.} {\it Dependence between Kinkelin and gamma-morphic functions}, Charkow Ges., 8 (1904),  123-135. JFM 34.0490.01

\bibitem[AlKu2023]{AK}
Alexanian S., Kuznetsov A.
{\it  On the Barnes double gamma function.}
Integral Transforms Spec. Funct. 34 (2023), no. 12, 891-914.

\bibitem[AAR1999]{AAR}
Andrews G. E., Askey R., Roy R.
{\it Special functions.} 
 Cambridge: Cambridge University Press., (1999).
 



\bibitem[App1881]{App}
 Appell P.
{\it Sur une classe de fonctions analogues aux fonctions Eul\'eriennes.} (French) 
Klein Ann. XIX, 84-102 (1881), JFM 13.0388.01.
	
\bibitem[Bakh1967]{Bah}
Bakhmutskaya \`E. Ya. {\it Mathematics in Kharkov university. Kharkov mathematical society.} (Russian) In {\it History of national mathematics in four volumes. Vol. 2, 1801-1917.} (eds. Z.~Shtokalo, A.~N.~Bogoljubov, 
Ju.~A.~Mitropol'ski\u\i, I.~B.~Pogrebysski\u\i, E.~Ja.~Remez, K.~A.~Rybnikov, Ju.~D.~Sokolov and V.~S.~Sologub), 460-472,
Naukova Dumka, Kiev, 1967.  
	
\bibitem[Bar1899]{Bar1}	
Barnes E. W. {\it The Genesis of the Double Gamma Functions}
Proc. Lond. Math. Soc. 31 (1899), 358-381. JFM 30.0389.03

\bibitem[Bar1900.1]{Bar2}
Barnes E. W.
{\it The theory of the $G$-function.} 
Quart. J. 31, 264-314 (1900). JFM 30.0389.02

\bibitem[Bar1900.2]{Bar3}
Barnes E. W.
{\it The theory of the double gamma Function.}  JFM 31.0441.01
Lond. R. S. Proc. 66, 265-266 (1900).

\bibitem[Bar1901]{Bar4}
Barnes E. W.
{\it The theory of the double gamma function.} 
Lond. Phil. Trans. (A) 196, 265-387 (1901). JFM 32.0442.02

\bibitem[Bar1903]{Bar-cap}
Barnes E. W.
{\it On the coefficients of capacity of two spheres.} 
Quart. J. 35, 155-175 (1903).  JFM 34.0913.01

\bibitem[Bar1904]{Bar5}
Barnes E. W. , {\it On the theory of the multiple gamma function}, Trans. Camb. Philos. Soc., 19 (1904), 374-425. JFM 35.0462.01

\bibitem[Bea1903]{Bea1903}
Beaupain J.
{\it Sur les fonctions d’ordre sup\'erieur de Kinkelin.} (French) JFM 34.0492.06
Belg. M\'em. C. et. sav. \'etr. 59, 67 S. (1903). JFM 34.0492.06



\bibitem[Bea1904]{Bea1904}
Beaupain J.
{\it Sur la fonction $\log\Gamma(a)$. Sur  la fonction $\log G_\lambda(a)$}. (French) 
13 S. 16 S. Belg. Bull. Science 62 (1904).
JFM 35.0461.03

\bibitem[Bea1907]{Bea1907}
Beaupain J.
{\it Sur la fonction gamma double.} (French) 
Belg. M\'em. et. sav. \'etr. in $4^\circ$
 (2) 1, 25 S. (1907). JFM 38.0469.04
 
\bibitem[Bel2000]{Bel}
Belomestnyh V. N., Belomestnyh L. A., {\it Physical-mathematical
education in high technical school in Siberia. Part 1}. (Russian) Tomsk Polytechnical Univ., 2000.



\bibitem[Ben1933]{Ben}
Bendersky L.
{\it Sur la fonction gamma g\'en\'eralis\'ee.} (French)
Acta Math. 61 (1933), 263-322 ,  JFM 59.0373.02.

\bibitem[Bou1960]{Bou}
Bourbaki N.
{\it \'El\'ements d'histoire des math\'ematiques.}
(French)
Histoire de la Pens\'ee, IV
Hermann, Paris, 1960. 

\bibitem[Erd1953]{Bei}
Erd\'elyi A., Magnus W., Oberhettinger F., Tricomi F. G.
{\it Higher transcendental functions. Vols. I, II} 
Based, in part, on notes left by Harry Bateman
McGraw-Hill Book Co., Inc., New York-Toronto-London, 1953.
Vol. III, 1955. 

\bibitem[Fad1995]{Fad}
Faddeev L. {\it Discrete Heisenberg-Weyl group and modular group}, Lett. Math. Phys. 34,
(1995), 249. 

\bibitem[Fock1997]{Fock0}
Fock V. V. {\it Dual Teichm\"uller spaces.}
Preprint\\
\vphantom{.}\hfill
{\tt https://arxiv.org/pdf/dg-ga/9702018v1.pdf}

\bibitem[FoCh1999]{Fock}
Fock V. V., Chekhov L. O.
{\it Quantum Teichm\"uller spaces.}
Theoret. and Math. Phys. 120 (1999), no.3, 1245-1259.


\bibitem[Gal2004]{Gal}
Galanova R. A.
{\it His thoughts and deeds were devoted to Siberia: about the first rector of TTI {\rm(}TPU{\rm)%
\footnote{ Tomsk Techological Institute (TTI)
changed couple of names, since 1991 it is  Tomsk Technological University,
(TTU).}} E.~L.~Zubashev.} (Russian). Tomsky polytehnik, 2004, 10, 160-163.

\bibitem[Gal2009]{Gal1}
Galanova R. A. {\it Efim Luk'yanovich Zubashev: a biographical essay.}
(Russian)
Tomsk Polytechnical University, Tomsk, 2009.


\bibitem[Gla1877]{Gla}
Glaisher, J. W. L.
{\it On the product $1^1\cdot 2^2\cdot 3^3\dots n^n$.}
Messenger (2) VII. 43-47 (1877). JFM 09.0190.01

\bibitem[Grad1963]{Gra1963}
Gradshteyn I. S., Ryzhik I. M.
{\it Table of integrals, series, and products.} (Russian)
Fourth edition prepared by Ju. V. Geronimus and M. Ju. Ce\u\i tlin.
Fizmatgiz, 1963. English transl.:
Academic Press, New York-London, 1965.

\bibitem[Gran2003]{Gra2003}
Granina Yu. G. {\it The professor-mathematician, the second director
of TTI V.~P.~Alekseevsky} (Russian), Izvestiya Tomskogo polytechnicheskogo universiteta,
[Bulletin of  Tomsk polytechnical university], 306, no. 3.,
(2003)  162-167.

\bibitem[Gus1950]{Gus0}
Gussov V. V.
{\it From the history of transcendental functions in Russia and USSR.}
(Russian)
Thesis (Ph.D.) -- Moscow State University, 1950, 385pp.




\bibitem[Gus1952]{Gus}
Gussov V. V.
{\it Works of Russian scholars on the theory of the gamma function.}(Russian) Istor.-Mat. Issled. 5 (1952), 421-472.

\bibitem[Gus1953]{Gus1}
Gussov, V. V.
{\it Development of the theory of cylinder functions in Russia and the USSR.} (Russian)
Istor.-Mat. Issled. (1953), 355-475.


\bibitem[Har1905.1]{Har1}
Hardy G. H.
{\it On the expression of the double zeta-function and double gamma-function in terms of elliptic functions.} 
Cambr. Phil. Trans. 20, 1-35 (1905). JFM 36.0508.03

\bibitem[Har1905.2]{Har2}
Hardy G. H.
{\it On double Fourier series, and especially those which represent the double zeta-function with real and incommensurable parameters.} 
Quart. J. 37, 53-79 (1905). JFM 36.0501.02.

\bibitem[Hei1847]{Hei}
 Heine E.
{\it Untersuchungen \"uber die Reihe 
$$1+\frac{(1-q^\alpha)(1-q^\beta)}{(1-q)(1-q^\gamma)}\cdot x+\frac{(1-q^\alpha)(1-q^{\alpha+1})(1-q^\beta)(1+q^{\beta+1})}
{(1-q)(1-q^2)}{(1-q^\gamma)(1+q^{\gamma+1})}\cdot x^2+\dots $$}
(German)
J. Reine Angew. Math. 34, 285-328 ({ 1847}),  ERAM%
\footnote{Electronic Research Archive for Mathematics.} 034.0971cj.

\bibitem[Hei1878]{Hei1}
Heine E.
{\it Handbuch der Kugelfunctionen. Theorie und Anwendungen.}
(German)
Berlin, G.~Reimer, 1878.  Reprinted Physica-Verlag, W\"urzburg, 1961.

\bibitem[Hol1886]{Hol}
H\"older O.
{\it Ueber eine transcendente Function.} (German) 
G\"ott. Nachr. 1886, 514-522 (1886). JFM 18.0376.01

\bibitem[KarKu2024]{KK}
Karp D.,  Kuznetsov A.
{\it Extending the Meijer G-function}. Preprint
\\
\vphantom{.}\hfill {\tt https://arxiv.org/abs/2403.05708.}

\bibitem[Kash1998]{Kash}
Kashaev R. M.
{\it Quantization of Teichm\"uller spaces and the quantum dilogarithm.}
 Math. Phys. 43 (1998), no.2, 105-115.
 
 \bibitem[KashN2011]{KN}
 Kashaev R. M., Nakanishi T.
{\it Classical and quantum dilogarithm identities.}
SIGMA Symmetry Integrability Geom. Methods Appl. 7 (2011), Paper 102, 29 pp.

 
 \bibitem[KLS2002]{KLS}
Kharchev S., Lebedev D., Semenov-Tian-Shansky M.
{\it Unitary representations of $U_q(\mathfrak{s}\mathfrak{l}(2,\mathbb{R}))$, the modular double and the multiparticle $q$-deformed Toda chains.}
Comm. Math. Phys. 225 (2002), no.3, 573-609.


\bibitem[Khar1904]{Khar}
 {\it Извлечение из протоколов заседаний} (Russian), Сообщ. Харьков. матем. общ. Вторая сер., 8 (1904), 283-287 [{\it Extracts from the minutes of meetings}, Charkow Ges., 8 (1904), 283-287 ]
 
 \bibitem[Kin1854]{Kin0}
Kinkelin H.
{\it Untersuchung \"uber die Formel}
$$
nF(nx)=f(x)+f\bigl(x+\tfrac1n\bigr)+f\bigl(x+\tfrac2n\bigr)+\dots 
f\bigl(x+\tfrac{n-1}n\bigr).
$$
(German),
Archiv der Mathematik und Physik, 22 (1854), 189-224. 

\bibitem[Kin1860]{Kin}
Kinkelin H.
{\it Ueber eine mit der Gammafunction verwandte Transcendente und deren Anwendung auf die Integralrechung.} (German)
J. Reine Angew. Math. 57, 122-138 (1860). ERAM 057.1509cj

\bibitem[Kru1967]{Kru}
Krulikovski\u\i\vphantom{.} N. N.
{\it History of the development of mathematics in Tomsk.} (Russian)
Izdat. Tomsk. Univ., Tomsk, 1967; Second edition 2006. 

\bibitem[Kru1967.1]{Kru1}
Krulikovski\u\i\vphantom{.}  N. N.
{\it Mathematics in  Tomsk Technological Institute.} (Russian)
In {\it History of national mathematics in four volumes. Vol. 2, 1801-1917.} (eds. Z.~Shtokalo, A.~N.~Bogoljubov, 
Ju.~A.~Mitropol'ski\u\i, I.~B.~Pogrebysski\u\i, E.~Ja.~Remez, K.~A.~Rybnikov, Ju.~D.~Sokolov and V.~S.~Sologub), 524-527.
Naukova Dumka, Kiev, 1967.  

\bibitem[Mari1978]{Mar}
Marichev O. I {\it The method of evaluation of integrals of special
functions {\rm(}theory and tables of formulas{\rm)}} (Russian),
Nauka i Tehnika, Minsk, 1978; English version:
Marichev O. I.
{\it Handbook of integral transforms of higher transcendental functions: theory and algorithmic tables.}   John Wiley, New York, 1983.

\bibitem[Mark1951]{Mark}
Markushevich A. I.
{\it Essays on the history of the theory of analytic functions}
(Russian).
Gosudarstv. Izdat. Tehn.-Teor. Lit., Moscow-Leningrad, 1951.

\bibitem[Mark1981]{Mark1}
Markushevich A. I. {\it 
Theory of analytic functions.} (Russian). In
{\it Mathematics of the XIXth century, Geometry, Theory of analytic functions} (eds. A.~N.~Kolmogorov, A.~P.~Yushkevich), Moscow 1981,
 115-269 (1981); English transl.: Laptev~B.~L., Rozenfeld~B.~A., Markushevich~A.~I.
{\it Mathematics of the 19th century.
Geometry, analytic function theory.} 
Birkh\"auser, Basel, 1996.

\bibitem[Mol1985]{Mol}
Molin F. \`E. [Molien Th.]
{\it Number systems.} (Russian)
Translated from the German by L.~A.~Bokut', N.~N.~Krulikovski\u\i\vphantom{.}
 and I.~V.~L'vov. Edited by A.~I.~Kostrikin. With appendices
  by Krulikovski\u\i, Bokut'\vphantom{.} and L'vov.
  Nauka, Sibirsk. Otdel., Novosibirsk, 1985. 126 pp.
  Zbl%
  \footnote{The database Zentralblatt f\"ur Mathematik.} 0607.01025

\bibitem[NiYo2009]{NY}
Nikeghbali A., Yor M. {\it The Barnes $G$-function and its relations with sums and products of generalized
Gamma convolution variables.} Electronic Communications in Probability. 2009; 14:396.

\bibitem[Ost2013]{Ost}
Ostrovsky D.
{\it Theory of Barnes beta distributions.}
Electron. Commun. Probab. 18, Paper No. 59, 16 p. (2013).

\bibitem[PBM1985]{PBM}
Prudnikov A. P., Brychkov Yu. A., Marichev O. I.
{\it Integrals and series.} Vol. 1-3: Gordon \& Breach 
1986-1990.

\bibitem[PBM1989]{PBM-plus}
Prudnikov A. P., Brychkov Yu. A., Marichev O. I.
{\it Calculation of integrals and the Mellin transform.}
(Russian) Mathematical analysis, Vol. 27, 3-146,
Progress in Science and Technology, VINITI, 1989;
Translated in J. Soviet Math. 54 (1991), no. 6, 1239-1341.


\bibitem[Raa1851]{Raa}
Raabe J. L. {\it Zur\"uckf\"uhrung einiger Summen und bestimmten Integrale auf die Jacob-Bernoullische Function.} (German)
J. Reine Angew. Math.,
1851(42), 348-367. ERAM 042.1172cj

\bibitem[Rim1947]{Rim1}
Rimski\u\i-Korsakov B. S.
{\it Notes on generalized multiplication theorems of Bernulli polynomials
and Kinkelin functions.} (Russian) Trudy MAI [Moscow Aviation Institute],  1, 1947


\bibitem[Rim1957]{Rim2}
Rimski\u\i-Korsakov B. S.
{\it A version of the construction of a theory of generalized gamma-functions based on the Laplace transform.}
Amer. Math. Soc. Transl. (2)17(1961), 201-217; Translated from: Moskov. Oblast. Ped. Inst. Uch. Zap. [Московский Областной Педагогический Институт. Ученые записки.] 57 (1957), 121-141.

\bibitem[Rim1957.1]{Rim3}
Rimski\u\i-Korsakov B. S.
{\it An investigation on theory of generalizations of the Gamma-function}
(Russian).
Dissertation  for the degree of candidate of physical and mathematical sciences,  Moscow, 1957.


\bibitem[Sar1987]{Sar}
Sarnak P.
{\it Determinants of Laplacians.} Comm. Math. Phys.110 (1987), no.1, 113-120.

\bibitem[Scho1884]{Sch}
Schobloch, J. A.
{\it \"Uber Beta- und Gammafunctionen.} (German)
Halle. Nebert (1884). JFM 16.0395.01.

\bibitem[Shi1976]{Shi1}
Shintani T.
{\it On Kronecker limit formula for real quadratic fields.} Proc. Japan Acad.52 (1976), no.7, 355-358.

\bibitem[Shi1977]{Shi2}
 Shintani T., {\it On a Kronecker limit formula for real quadratic fields}, J. Fac. Sci. Univ.
Tokyo, Sect. 1A 24 (1977), 167-199

\bibitem[Sin1917]{Sin}
Синцов Д. М. {\it В. П. Алексеевский (1858-1916)}. (Russian) Сообщ. Харьков. матем. общ. Вторая сер., 15:5-6 (1917),  288-295
[Sintsov D. M. {\it V. P. Alexeevsky (1858-1916).}
Charkow Ges., 15, no. 5-6 (1917),  288-295].

\bibitem[Srin2007]{Sri}
Srinivasan G. K.
{\it The gamma function: an eclectic tour.}
Amer. Math. Monthly 114 (2007), no. 4, 297-315.	 
	 

\bibitem[Ste1904]{Ste}
 Stekloff  W. A. [Steklov V. A.] {\it Remarques relatives aux formules sommatoires d'Euler et de Bool} (French), Сообщ. Харьков. матем. общ. Вторая сер. [Charkow Ges.], 8 (1904), 136-195. JFM 35.0263.01
 
 \bibitem[Ste1991]{Ste91} 
 Steklov V. A. {\it Correspondence with national mathematicians. Recollections.} (Russian), eds.  V.~S.~Vladimirov, E.~P.~Ozhigova, V.~S.~Sobolev, Leningrad, Nauka, 1991.
	 
\bibitem[StSi1904]{StSi}
Стеклов В. А., Синцов Д. М. 
{\it  Отзыв профессоров Стеклова и Синцова об ученых трудах прив.-доцента
В.~П.~Алексеевского.} (Russian) Ученые записки Императорского Харьковского 
университета, 1904, 1, 1-11
[{\it Review of professors Steklov and Sintsov on the scientific works of privatdozent
V.~P.~Alekseevsky.} Scientific notes of the Imperator's Kharkov University,
1904, 1, 1-11]

\bibitem[Var1988]{Var}
Vardi I.  {\it Determinants of Laplacians and Multiple Gamma Functions}. SIAM Journal on Mathematical Analysis, 19(2) (1988), 493-507.
	 
	 \bibitem[Vig1979]{Vig}
 Vign\'eras  M. F. {\it  L'\'equation fonctionelle de la fonction z\^{e}ta de Selberg du groupe modulaire $SL(2,\Z)$} (French), Ast\'erisque 61 (1979), 235-249 	


 
 \bibitem[Wae1985]{Wae}
van der Waerden B. L.
{\it A history of algebra.
From al-Khw\=arizm\=\i\vphantom{.} to Emmy Noether.}
Springer-Verlag, Berlin, 1985. 
 
 \bibitem[Whi1954]{Whi}
 Whittaker E. T.
{\it Ernest William Barnes, 1874-1953}. Obit. Notices Roy. Soc. London 9 (1954), 15-25.
 
 \bibitem[WhWa1912]{WhWa}
 Whittaker  E. T.,  Watson G. N., {\it A course of Modern Analysis}, Second edition, 1912.

\bibitem[Wor2000]{Wor}
Woronowicz S.L. {\it Quantum exponential function}, Rev. Math. Phys. 12 (2000), 873-920.


 \bibitem[Yush1968]{Yush}
Yushkevich A. P.
{\it History of mathematics in Russia up to 1917.} (Russian)
Nauka, Moscow, 1968, 591~pp.

\end{thebibliography}
\end{document}